\documentclass[12pt]{article}

\usepackage[utf8]{inputenc}
\usepackage[T1]{fontenc} 

\usepackage{authblk}
\usepackage{amsmath,amsthm,amssymb,amsfonts}
\usepackage{mathtools}
\usepackage{microtype}
\usepackage{enumitem}
\usepackage{multirow,multicol}
\usepackage{colortbl}
\usepackage{caption}
\usepackage{graphicx}
\usepackage{longtable}

\usepackage{todonotes}
\usepackage{adjustbox}
\usepackage{float}
\usepackage{tikz}
\usetikzlibrary{arrows.meta}
\usepackage{pdfpages}

\usepackage{hyperref} 
\hypersetup{
   pdftitle={Unramified abelian covers with many points},    
   pdfauthor={Jean Gasnier},            
   pdfsubject={},           
   colorlinks=true,         
   linkcolor=teal,         
   citecolor=purple,         
}
\usepackage{url}
\usepackage[onelanguage,tworuled,noline,noend,algosection,linesnumbered,french]{algorithm2e}

\numberwithin{equation}{section}

\newtheorem{theo}{Theorem}

\theoremstyle{definition}

\newtheorem{prop}[theo]{Proposition}

\theoremstyle{remark}

\definecolor{vertpomme}{RGB}{130,210,30}

\newcommand{\F}{\ensuremath{\mathbb{F}}}

\newcommand{\J}{\mathcal{J}}

\DeclareMathOperator{\Id}{Id} 

\title{Unramified abelian covers with many points}
\date{\today}
\author{Jean Gasnier}
\affil{Univ. Bordeaux, INRIA, CNRS, Bordeaux INP, INRIA, IMB, UMR 5251,  F-33400, Talence,France}

\begin{document}
\maketitle

\begin{abstract}
We produce curves with a record number of points over the finite fields with $4$, $9$, $16$ and $25$ elements, as unramified abelian covers of algebraic curves.
\end{abstract}

\section{Introduction}

The observation that the Weil bound stops being tight as the genus grows over a fixed finite field has motivated further investigations of the maximal number of points a smooth projective curve of given genus can have. Other than the intrinsic interest of this question, this value is key in the theory of algebraic-geometry codes, as it determines the quality of such codes (see for example~\cite{Sti08}).

One can take several approaches to this question. A first approach is to improve the Weil bound, and obtain a formula providing a better estimate for larger genera. Famous results in this direction were found by Oesterlé (see~\cite[Chapter VI]{Serre20}) and Serre (see~\cite[Chapter II]{Serre20}). A second approach consists in giving an estimate of the Ihara constant
\[A(q) = \limsup_{X/\F_q}\frac{\#X(\F_q)}{g_X}\]
where $g_X$ denotes the genus of the curve $X$. The Ihara constant is known to be $\sqrt{q}-1$ when $q$ is a square thanks to a contribution by Drinfeld and Vladut~\cite{VlaDri83}, which shows that $A(q)\leqslant \sqrt{q}-1$ for any prime power $q$, and a contribution by Ihara~\cite{Ihara81} which provides an example of a family of curves reaching this bound when $q$ is a square.
The third approach consists in computing this value for a fixed genus over a fixed finite field. This requires computing actual examples of curves with many points. Serre's course~\cite{Serre20} gives a summary of most known methods used to produce such curves. In particular, in Section~5.12 for instance, he mentions how one can obtain curves with many points as abelian covers of curves with lower genus, using Class Field Theory (CFT).

The manYPoints database~\cite{manypoints} provides a compilation of results on this matter for small genera ($g\leqslant 50$) over small finite fields. For a given pair $(q,g)$, it provides the largest known number of points posessed by a curve of genus $g$ over $\F_q$, together with the smallest known upper bound on this value. 

In this paper, we use CFT in the unramified case to produce new curves with many points, improving previous records from~\cite{manypoints}. In Section~\ref{sec-recalls-cft}, we make some recalls on CFT in the unramified case. In Section~\ref{sec-record}, we explain how we used the data from LMFDB~\cite{lmfdb} to compute our new curves. 

\section{Some recalls on abelian covers}\label{sec-recalls-cft}

Let $K$ be a finite field with $q$ elements, and let $X$ be a smooth projective curve over $K$. Note that in this article, every curve is assumed to be \textbf{geometrically integral}.
Let $P$ be a rational point of $X$. We intend to describe all the unramified abelian covers of $X$ totally split above $P$, up to isomorphism. 
This can be achieved using Class Field Theory (CFT), where one uses the Jacobian variety $\J_X$ of $X$ to describe these covers. 

Let $Y$ be a smooth projective curve over $K$. One says that $Y$ is a cover of $X$ if there exists a non-zero finite morphism 
\[\tau : Y \longrightarrow X.\]
Note that $\tau$ induces a non-zero morphism of $K$-algebras $ K(X)\longrightarrow K(Y)$, by~\cite[Section 7, Proposition 3.13]{Liu02}.
Thus, $K(Y)$ is a finite extension of $K(X)$. One says that $Y$ is an unramified abelian cover of $X$ if $K(Y)/K(X)$ is an unramified abelian extension. In this context, we call $\tau$ an unramified abelian covering map of $X$.

The main result that we are going to use is that any such cover is the pullback of a separable isogeny $\theta$ by a regular map from $X$ to $\J_X$. One may refer to~\cite{Serre84} for further details, where this result is stated in Chapter~1. Since we are only interested in unramified abelian covers that are totally split above $P$, it is sufficient to consider pullbacks of isogenies by $j_P$, the Jacobi map associated to $P$, i.e.\ the closed immersion of $X$ in $\J_X$ mapping $P$ to $0$. In this context, it is equivalent to ask that $P$ splits completely in $Y$, and that $\ker \theta$ is entirely composed of rational points. One then infers Theorem~\ref{theo-hilbert-function-geom}.
\begin{theo}[\protect{\cite[Chapter I, Corollary of Theorem 4 and Theorem 5]{Serre84}}]\label{theo-hilbert-function-geom}
Let $K$ be a finite field and let $X$ be a smooth projective curve over $K$. Let $P$ be a rational point of $X$.
There exists a smooth projective curve $Y_{max}$ over $K$ and an abelian covering 
\[\tau_{max} : Y_{max}\longrightarrow X\]
of $X$ that is unramified, totally split above $P$, and maximal in the following sense:
for any unramified abelian covering
\[\tau: Y\longrightarrow X\] 
satisfying these properties, there exists an unramified abelian covering 
\[\tau_{Y_{max}/Y}:Y_{max}\longrightarrow Y\] 
such that 
\[\tau_{max} = \tau\circ\tau_{Y_{max}/Y}.\]
The cover $Y_{max}$ is obtained by pulling back the isogeny
\[ \varphi = F_{\J_X}-\Id,\]
by the Jacobi map $j_P$. Its Galois group is isomorphic to $\J_X(K)$.
\end{theo}

Theorem~\ref{theo-hilbert-function-geom} induces a correspondence between the isomorphism classes of unramified abelian covers of $X$ totally split above $P$ and the subgroups of $\J_X(K)$. Indeed, let $H$ be a subgroup of $\J_X(K)$. Let 
\[\pi_H : \J_X \longrightarrow \J_X/H\]
be the quotient isogeny associated to $H$.
Then there exists \[\theta :\J_X/H\longrightarrow\J_X \] an isogeny such that \[\theta \circ \pi_H = F_{\J_X}-\Id.\]
The cover obtained by pulling back $\theta$ by $j_P$ is an unramified abelian cover of $X$, totally split above $P$, with Galois group $\J_X(K)/H$.
Figure~\ref{fig-cft} sums up the situation.
\begin{figure}[H]
\centering
\begin{tikzpicture}
\node (Ymax) at (0,4) {$Y_{max}$};
\node (Y) at (0,2) {$Y$};
\node (X) at (0,0) {$X$};
\node (Jx) at (2,0) {$\J_X$};
\node (JxH) at (2,2) {$\J_X/H$};
\node (Jxm) at (2,4) {$\J_X$};
\node[below] (jp) at (1,0) {$j_P$};
\node[left] (tau) at (0,1) {$\tau$};
\node[left] (theta) at (2,1) {$\theta$};
\node[left] (piH) at (2,3) {$\pi_H$};

\draw[->] (Ymax)--(Y);
\draw[->] (Y)--(X);
\draw[->] (X)--(Jx);
\draw[->] (Y)--(JxH);
\draw[->] (Ymax)--(Jxm);
\draw[->] (Jxm)--(JxH);
\draw[->] (JxH)--(Jx);
\draw[->] (Jxm) to[bend left=45] node[midway,right]{$F_{\J_X}-\Id$} (Jx);
\end{tikzpicture}
\caption{Unramified abelian covers of $X$ totally split above $P$}\label{fig-cft}
\end{figure}

Now, Proposition~\ref{prop-decompo-geom} describes which rational points of $X$ are totally split in $Y$, the cover associated to $H\subset \J_X(K)$.
\begin{prop}\label{prop-decompo-geom}
Using the notation from Theorem~\ref{theo-hilbert-function-geom}, let $Q$ be a rational point of $X$, and let $H$ be a subgroup of $\J_X(K)$. Let $\tau: Y \longrightarrow X$ be the unramified abelian cover, totally split above $P$, associated with $H$.
Then $Y$ is totally split above $Q$ if and only if $j_P(Q)\in H$.
\end{prop}
Indeed $Q$ is totally split in $Y$ if and only if the fiber of $\theta$ above $j_P(Q)$ is entirely composed of rational points, which is only the case when $j_P(Q)\in H$.

\section{New record curves}\label{sec-record}

Let $K$ be a finite field with $q$ elements and let $X$ be a smooth projective curve over $K$. Assume there exists $P$ a rational point of $X$. Our goal is to construct a smooth projective curve $Y$ over $K$, with many rational points relatively to its genus, as an unramified abelian cover of $X$ totally split above $P$. Therefore, let $H\subset \J_X(K)$ be a subgroup of rational points of the Jacobian of $X$, and let $Y$ be the associated cover of $X$ as defined in Figure~\ref{fig-cft}, with Galois group isomorphic to $\J_X(K)/H$.

Let $g_X$ be the genus of $X$. Then the Riemann-Hurwitz formula allows us to express the genus $g_Y$ of $Y$:
\[g_Y-1 = \frac{|\J_X(K)|}{|H|}(g_X-1).\]
Moreover, any rational point $Q_Y$ of $Y$ must be an element of the fiber of a rational point $Q$ of $X$ that totally splits in $Y$. Any such fiber is composed of $|\J_X(K)/H|$ elements. Thus, by Proposition~\ref{prop-decompo-geom},
\[\#Y(K) =  \frac{|\J_X(K)|}{|H|}\#\{Q \in X(K)| j_P(Q)\in H\}.\]

It is a difficult task to find appropriate pairs $(X,H)$ defining curves $Y$ with many points. In what follows, we will rely on a trick, previously exploited in~\cite{GurXin22,NieXin98,Quebbemann89,vanderGeer09}. Assume that $K$ is an extension of degree 2 of a field $\kappa$ with $\sqrt{q}$ elements. Let $X_\kappa$ be a smooth projective curve over $\kappa$, and let $P$ be a $\kappa$-rational point of $X_\kappa$. Take
\begin{itemize}
    \item $X=(X_\kappa)_K$.
    \item $H = \J_{X_\kappa}(\kappa) \subset \J_X(K)$.
\end{itemize}
One can see that in this context, the subset of rational points of $X$ totally split in $Y$ is $X_\kappa(\kappa)$. Thus the data of the L-polynomial of $X_\kappa$ alone allows us to compute $\#X_\kappa(\kappa)$, the order $|\J_{X_\kappa(\kappa)}|$, and $|\J_X(K)|$, and therefore $g_Y$ and $\#Y(K)$.

LMFDB~\cite{lmfdb} lists, among other things, L-polynomials of algebraic curves over finite fields with $2$, $3$, $4$, and $5$ elements (among others). By enumerating this data, we can produce new curves with a record number of points over finite fields with $4$, $9$, $16$, and $25$ elements. These records are presented in Tables~\ref{tab-record-F4},~\ref{tab-record-F9},~\ref{tab-record-F16}, and~\ref{tab-record-F25}.
\begin{table}[htbp]
  \[
  \begin{array}{|c|c|c|c|c|}
    \hline
    \text{\'LMFDB label of }X & |G|  &  g_{Y} & \#Y(\F_4) & \text{Old record}~(\cite{manypoints}) \\
    \hline
    \text{\href{https://www.lmfdb.org/Variety/Abelian/Fq/4/2/d_i_o_x}{4.2.d\_i\_o\_x}} & 11 & 34 & 66 & 65\\
    \text{\href{https://www.lmfdb.org/Variety/Abelian/Fq/5/2/e_m_ba_bv_cu}{5.2.e\_m\_ba\_bv\_cu}} & 12 & 49 & 84& 81 \\
    \hline
  \end{array}
  \]
\caption{New curves with a record number of rational points over $\F_4$}\label{tab-record-F4}
\end{table}

\begin{table}[htbp]
  \[
  \begin{array}{|c|c|c|c|c|}
    \hline
    \text{\'LMFDB label of }X & |G|  &  g_{Y} & \#Y(\F_9) & \text{Old record}~(\cite{manypoints}) \\
    \hline
    \text{\href{https://www.lmfdb.org/Variety/Abelian/Fq/4/3/i_bi_ds_hn}{4.3.i\_bi\_ds\_hn}} & 9 & 28 & 108 & 105\\
    \text{\href{https://www.lmfdb.org/Variety/Abelian/Fq/4/3/h_ba_co_ez}{4.3.h\_ba\_co\_ez}} & 11 & 34 & 121 & 114\\
    \text{\href{https://www.lmfdb.org/Variety/Abelian/Fq/4/3/h_bb_ct_fk}{4.3.h\_bb\_ct\_fk}} & 12 & 37 & 132 & 126\\
    \hline
  \end{array}
  \]
\caption{New curves with a record number of rational points over $\F_9$}\label{tab-record-F9}
\end{table}

\begin{table}[htbp]
  \[
  \begin{array}{|c|c|c|c|c|}
    \hline
    \text{\'LMFDB label of }X & |G|  &  g_{Y} & \#Y(\F_{16}) & \text{Old record}~(\cite{manypoints})\\
    \hline
    \text{\href{https://www.lmfdb.org/Variety/Abelian/Fq/3/4/g_v_bx}{3.4.g\_v\_bx}} & 19 & 39 & 209 & 194\\
    \text{\href{https://www.lmfdb.org/Variety/Abelian/Fq/3/4/f_p_bg}{3.4.f\_p\_bg}} & 23 & 47 & 230 & \emptyset\\
    \hline
  \end{array}
  \]
\caption{New curves with a record number of rational points over $\F_{16}$}\label{tab-record-F16}
\end{table}

\begin{table}[htbp]
  \[
  \begin{array}{|c|c|c|c|c|}
    \hline
    \text{\'LMFDB label of }X & |G|  &  g_{Y} & \#Y(\F_{25}) & \text{Old record}~(\cite{manypoints}) \\
    \hline
    \text{\href{https://www.lmfdb.org/Variety/Abelian/Fq/3/5/k_bv_fc}{3.5.k\_bv\_fc}} & 16 & 33 & 256 & 226\\
    \text{\href{https://www.lmfdb.org/Variety/Abelian/Fq/3/5/j_bn_ec}{3.5.j\_bn\_ec}} & 20 & 41 & 300 & 260 \\
    \text{\href{https://www.lmfdb.org/Variety/Abelian/Fq/3/5/j_bo_eh}{3.5.j\_bo\_eh}} & 21 & 43 & 315 & 276\\
    \text{\href{https://www.lmfdb.org/Variety/Abelian/Fq/3/5/i_bf_dc}{3.5.i\_bf\_dc}} & 24 & 49 & 336 & 315\\
    \hline
  \end{array}
  \]
\caption{New curves with a record number of rational points over $\F_{25}$}\label{tab-record-F25}
\end{table}

\newpage
\bibliographystyle{alpha}
\bibliography{biblio.bib}
\end{document}